\documentclass[12pt]{article} 
\usepackage[utf8]{inputenc} 
\usepackage{geometry}
\usepackage{graphicx} 
\usepackage{amssymb}
\usepackage{amsmath}
\usepackage{cite}
\usepackage{mathrsfs}
\usepackage{booktabs} 
\usepackage{array} 
\usepackage{paralist}
\usepackage{verbatim} 
\usepackage{subfig}
\usepackage{fancyhdr} 
\usepackage{sectsty}
\allsectionsfont{\sffamily\mdseries\upshape} 
\usepackage[nottoc,notlof,notlot]{tocbibind} 
\usepackage[titles,subfigure]{tocloft} 

\title{Polynomial solutions of the Dirichlet problem for the Tricomi equation in a strip.}
\author{\bf Oleg D. Algazin}
\date{Bauman Moscow State Technical University, Moscow, Russia} 

\begin{document}
\maketitle
\thispagestyle{empty}

\begin{abstract}
         An inhomogeneous Tricomi equation is considered in a strip with a polynomial right-hand side. It is shown that the Dirichlet boundary value problem with polynomial boundary conditions has a polynomial solution. An algorithm for constructing this polynomial solution is given and examples are considered. If the strip lies in the ellipticity region of the equation, then this solution is unique in the class of functions of polynomial growth. If the strip lies in a mixed region, then the solution of the Dirichlet problem is not unique in the class of functions of polynomial growth, but it is unique in the class of polynomials.
\end{abstract}


\section*{Introduction}

The equation
\[
yu_{xx}+u_{yy}=0
\]
 was first considered by F.Tricomi \cite{Tr} and subsequently received his name. This equation is elliptic in the upper halfplane $y> 0$, hyperbolic in the lower halfplane $y <0$, and parabolically degenerates on the straight line $y = 0$. Equations of a mixed type are used in transonic gas dynamics \cite{B}-\cite{B1}. The Dirichlet problem for a mixed-type equation in a mixed domain is, in general, ill-posed  \cite{B2}. A lot of works has been devoted to finding conditions for the well-posedness  of the formulation of the Dirichlet problem for a mixed-type equation in a mixed domain, for example \cite{N}-\cite{Sol}.

      This paper is devoted to finding exact polynomial solutions of the inhomogeneous Tricomi equation in a strip with a polynomial right-hand side. By means the Fourier transform, in the same way as was done in \cite{A} for the Poisson equation, it is shown that the Dirichlet boundary value problem with polynomial boundary conditions has a polynomial solution. An algorithm for constructing this polynomial solution is given and examples are considered. If the strip lies in the ellipticity region of the equation, then this solution is unique in the class of functions of polynomial growth. If the strip lies in a mixed domain, then the solution of the Dirichlet problem is not unique in the class of functions of polynomial growth, but it is unique in the class of polynomials.

\section{Statement of the problem.  }

       We consider the Tricomi equation in a strip with a polynomial right-hand side
\[
                  Tu(x,y):=yu_{xx} (x,y)+u_{yy} (x,y)=P(x,y),~~   a<y<b, \eqno{(1)}               
\]
where $ P(x,y)$ is a polynomial in the variables $x$ and $y$.

     On the boundary of the strip we give the Dirichlet boundary conditions
\[
                         u(x,a)=\phi(x),~~u(x,b)=\psi(x), \eqno{(2)}
\]                                                      
where $\phi(x)$  and $\psi(x)$ are polynomials.

     If $\tilde u(x,y)$ is a particular polynomial solution of the inhomogeneous Tricomi equation (1), then for the function $v(x,y)=u(x,y)-\tilde u(x,y)$, we obtain the homogeneous equation
\[
                              Tv(x,y)=0,~~    a<y<b,   \eqno{(3)}
\]
and the Dirichlet boundary conditions
\[
                    v(x,a)=\phi(x)-\tilde u(x,a),~~v(x,b)=\psi(x)-\tilde u(x,b).  \eqno{(4)}
\]
Solving the Dirichlet problem for the homogeneous Tricomi equation (3),(4), we obtain solution of the Dirichlet problem for the inhomogeneous Tricomi equation (1),(2) according to the formula
\[
u(x,y)=v(x,y)+\tilde u(x,y).
\]

      The solution of the Dirichlet problem will be sought in the class of functions of polynomial growth with respect to $x$:
\[
\int_{\mathbb{R}}|u(x,y)|(1+|x|)^{-m}dx<C  \eqno{(5)}
\]
                                           
for some $m\geq 0$ and for each $y\in(a,b)$.

\section{A particular solution of the Tricomi equation}

The Tricomi equation with the polynomial right-hand side $P(x,y)$,
\[
Tu(x,y):=yu_{xx} (x,y)+u_{yy} (x,y)=P(x,y),  
\]
has polynomial solutions, one of which can be obtained from the formula given below. It is sufficient to give this formula for a monomial $P(x,y)=x^n y^m$.  Then a particular solution of the Tricomi equation with the right-hand side $x^n y^m$ is the function
\[
\tilde u(x,y)=T^{-1} (x^n y^m )=
\]
\[ =\sum_{j=0}^{[n/2]}(-1)^j\frac{n!(m-1)!\prod_{k=0}^j(m+3k)}{(m+3j+2)!(n-2j)!} y^{m+3j+2} x^{n-2j}  ,       \eqno{(6)} 
\]           
where  $[n/2]$ is the integer part of the number $n/2$.

     The validity of formula (6) is proved by direct verification
\[
T\tilde u(x,y)=x^n y^m.
\]
For example, for
\[
P(x,y)= x^5 y^7
\]
a particular solution according to (6) will be
\[
\tilde u(x,y)=\frac{1}{72} x^5 y^9-\frac{5}{2376} x^3 y^{12}+\frac{1}{16632} xy^{15}.
\]

\section{The Dirichlet problem for a strip lying in the domain of ellipticity}

       We consider the strip $0<y<b$, in which the Tricomi equation is elliptic and parabolically degenerates on the boundary line $y=0$. Since the solution of the Dirichlet problem for an inhomogeneous equation reduces to solving the Dirichlet problem for a homogeneous equation, we consider the Dirichlet problem for a homogeneous equation.
\[
                             Tu(x,y):=yu_{xx} (x,y)+u_{yy} (x,y)=0,~~0<y<b,    \eqno{(7)}
\]
\[
               u(x,0)=\phi(x),~~   u(x,b)=\psi(x),                                        \eqno{ (8)}
\]
where $\phi(x)$ and $\psi(x)$ are polynomials.

     Since the functions $u(x,y)$  of slow growth with respect to $x$ determine regular functionals from the space  of generalized functions of slow growth $\mathscr S^{'}  (\mathbb R  )$ for each $y\in(0,b)$, we can apply the Fourier transform with respect to $x$\cite{V}:
  \[   
\mathscr F_x [u(x,y) ](t,y)=U(t,y).
\]
We apply the Fourier transform with respect to $x$ to equation (7) and the boundary conditions (8). We obtain a boundary value problem for a second-order ordinary differential equation with parameter $t\in\mathbb R$ 
\[
        -t^2 yU(t,y)+U_{yy} (t,y)=0,~~0<y<b       \eqno{(9)}      
\]                         
and boundary conditions
\[
       U(t,0)=\Phi(t),~~U(t,b)=\Psi(t).            \eqno{(10)}                                 
\] 

Equation (9) is the Airy equation, its general solution is expressed in terms of Airy functions
\[
U(t,y)=c_1 (t)\textup{Ai}(t^{2/3} y)+c_2 (t)\textup{Bi}(t^{2/3} y).
\]
Using the boundary conditions (10), we obtain a unique solution of the boundary value problem (9), (10)
\[
          U(t,y)=K(t,y)\Phi(t)+L(t,y)\Psi(t),      \eqno{(11)}                                     
\]
where
\[
K(t,y)=\frac{\textup{Ai}(t^{2/3} y)\textup{Bi}(t^{2/3} b)-\textup{Bi}(t^{2/3} y)\textup{Ai}(t^{2/3} b)}{ \textup{Ai}(0)\textup{Bi}(t^{2/3} b)-\textup{Bi}(0)\textup{Ai}(t^{2/3} b) },
\]
\[
L(t,y)=\frac{\textup{Ai}(0)\textup{Bi}(t^{2/3} y)-\textup{Bi}(0)\textup{Ai}(t^{2/3} y)}{ \textup{Ai}(0)\textup{Bi}(t^{2/3} b)-\textup{Bi}(0)\textup{Ai}(t^{2/3} b)}.
\]
Applying the inverse Fourier transform, we obtain a unique solution of the Dirichlet problem (7), (8) in the class of functions of polynomial growth in the form of convolution
\[
       u(x,y)=k(x,y)*\phi(x)+l(x,y)*\psi(x),   \eqno{(12)}                                     
\]
where
\[
k(x,y)=\mathscr F_t^{-1} (K(t,y) ), ~~  l(x,y)=\mathscr F_t^{-1} (L(t,y) ).
\]

To find the convolution (12) with the polynomials $\phi(x)$ and $\psi(x)$, it suffices to consider the case of a monomial $x^n$.

$K(t,y)$ and $L(t,y)$ are infinitely differentiable and rapidly decreasing as functions of the variable $t$, that is, they belong to the space $\mathscr S(\mathbb R)$, and hence both $k(x,y)$ and $l(x,y)$ as functions of the variable $x$ belong to the space $\mathscr S(\mathbb R)$, since the Fourier  transform  translates $\mathscr S(\mathbb R)$ into itself. In addition, $K(t,y)$ and $L(t,y)$ are even functions of the variable $t$ and, therefore, $k(x,y)$ and $l(x,y)$ are even functions of the variable $x$.

The solution of the Dirichlet problem (7), (8) with $\phi(x)=x^n,\psi(x)=0$  is the function
\begin{gather*}
u_n(x,y)=k(x,y)*x^n=\int_{-\infty}^{\infty}(x-t)^nk(t,y)dt=\\
=\int_{-\infty}^{\infty}\sum_{j=0}^nC_n^jx^{n-j}t^j(-1)^jk(t,y)dt=\sum_{j=0}^nC_n^jx^{n-j}(-1)^j\int_{-\infty}^{\infty}t^jk(t,y)dt,
\end{gather*}
where $C_n^j=n!/j!(n-j)!$ are binomial coefficients. Since the last integral for odd $j$ is zero because $k(t,y)$ is even in $t$, then
\[
u_n (x,y)=\sum_{m=0}^{[n/2]}C_n^{2m} x^{n-2m} \int_{-\infty}^{\infty}t^{2m}k(t,y)dt,
\]
where  $[n/2]$ is the integer part of the number $n/2$. Using the properties of the Fourier transform, we obtain
\begin{gather*}
p_{2m}(y)=\int_{-\infty}^{\infty}t^{2m}k(t,y)dt=\lim_{x\to 0}\int_{-\infty}^{\infty}t^{2m}k(t,y)e^{ixt}dt=\\
=\lim_{x\to 0}\mathscr{F}_t\left[t^{2m}P_1(t,y)\right](x,y)=(-1)^m\lim_{x\to 0}\frac{d^{2m}}{dx^{2m}}K(x,y).
\end{gather*}
$p_{2m}(y)$ are polynomials in $y$ whose generating function is $K(t,y)$.
\[
K(t,y)=\frac{\textup{Ai}(t^{2/3} y)\textup{Bi}(t^{2/3} b)-\textup{Bi}(t^{2/3} y)\textup{Ai}(t^{2/3} b)}{ \textup{Ai}(0)\textup{Bi}(t^{2/3} b)-\textup{Bi}(0)\textup{Ai}(t^{2/3} b) }=\sum_{m=0}^{\infty}p_{2m}(y)\frac{(-1)^m t^{2m}}{(2m)!}.
\]

We give some first polynomials $p_{2m}(y)$   and the corresponding solutions of the Dirichlet problem $u_n (x,y)$.
\begin{gather*}
p_0 (y)=\frac{b-y}{b},\\
p_2 (y)=\frac{y}{6b} (y^3-2by^2+b^3 ),\\
p_4 (y)=\frac{y}{210b} (-10y^6+28by^5-35b^3 y^3+17b^6 ),\\
p_6 (y)=\frac{y}{252b} (4y^9-14by^8+30b^3 y^6-51b^6 y^3+31b^9 )\\
p_8 (y)=\frac{y}{54054b} (-308y^{12}+1274by^{11}-4004b^3 y^9+14586b^6 y^6-31031b^9 y^3+19483b^{12} ),\\
p_{10} (y)=\frac{y}{36036b} (77y^{15}-364by^{14}+1540b^3 y^{12}-9724b^6 y^9+\\+44330b^9 y^6-97415b^{12} y^3+61556b^{15} ).\\
u_0 (x,y)=\frac{b-y}{b},\\
u_1 (x,y)=\frac{x(b-y)}{b},\\
u_2 (x,y)=\frac{1}{6b} (-6x^2 y+6x^2 b+y^4-2by^3+b^3 y),\\
u_3 (x,y)=\frac{x}{2b} (-2x^2 y+2x^2 b+y^4-2by^3+b^3 y),\\
u_4 (x,y)=\frac{1}{210b}(-210x^4 y+210x^4 b+210x^2 y^4-420x^2 y^3b+\\+210x^2 yb^3-10y^7+28by^6-35b^3 y^4+17b^6 y),\\
u_5 (x,y)=\frac{x}{42b} (-42x^4 y+42x^4 b+70x^2 y^4-140x^2 y^3 b+70x^2 yb^3-\\-10y^7+28by^6-35b^3 y^4+17b^6 y).
\end{gather*}

Similarly, the solution of the Dirichlet problem (7), (8) with $\phi(x)=0,\psi(x)=x^n$ is the function
\[   
v_n (x,y)=\sum_{m=0}^{[n/2]}C_n^{2m} x^{n-2m}q_{2m}(y),
\]
where $q_{2m}(y)$   are polynomials in $y$ whose generating function is $L(t,y)$.
\[
L(t,y)=\frac{\textup{Ai}(0)\textup{Bi}(t^{2/3} y)-\textup{Bi}(0)\textup{Ai}(t^{2/3} y)}{ \textup{Ai}(0)\textup{Bi}(t^{2/3} b)-\textup{Bi}(0)\textup{Ai}(t^{2/3} b)}=\sum_{m=0}^{\infty}q_{2m}(y)\frac{(-1)^m t^{2m}}{(2m)!}.
\]

We give some first polynomials $q_{2m}(y)$   and the corresponding solutions of the Dirichlet problem $v_n (x,y)$.
\begin{gather*}
 q_0 (y)=\frac{y}{b},\\
q_2 (y)=\frac{y}{6b} (b^3-y^3 ),\\
q_4 (y)=\frac{y}{42b} (2y^6-7b^3 y^3+5b^6 ),\\
q_6 (y)=\frac{y}{252b} (-4y^9+30b^3 y^6-75b^6 y^3+49b^9 ),\\
q_8 (y)=\frac{y}{4914b} (28y^{12}-364b^3 y^9+1950b^6 y^6-4459b^9 y^3+2845b^{12} ),\\
q_{10} (y)=\frac{y}{3276b} (-7y^{15}+140b^3 y^{12}-1300b^6 y^9+6370b^9 y^6-14225b^{12} y^3+9022b^{15} ).\\
v_0 (x,y)=\frac{y}{b},\\
v_1 (x,y)=\frac{xy}{b},\\
v_2 (x,y)=\frac{y}{6b} (6x^2-y^3+b^3 ),\\
v_3 (x,y)=\frac{xy}{2b} (2x^2-y^3+b^3 ), \\
v_4 (x,y)=\frac{y}{42b} (42x^4-42x^2 y^3+42x^2 b^3+2y^6-7b^3 y^3+5b^6 ),\\
v_5 (x,y)=\frac{xy}{42b} (42x^4-70x^2 y^3+70x^2 b^3+10y^6-35b^3 y^3+25b^6 ).
\end{gather*}

\textbf{Example . } We consider the Dirichlet problem for the inhomogeneous Tricomy equation
\[
Tu(x,y)=2x^2 y^3-18x^3 y^4,~~x\in\mathbb R~~ ,0<y<b,
\]
\[
u(x,0)=0,~~u(x,b)=0,~~x\in\mathbb R.
\]
A particular solution of the equation is the polynomial
\[
\tilde u(x,y)=-\frac{3}{5} x^3 y^6+\frac{1}{10} x^2 y^5+\frac{1}{20} xy^9-\frac{1}{280} y^8
\]
and the Dirichlet problem for an inhomogeneous equation reduces to the Dirichlet problem for a homogeneous equation for the function $v(x,y)=u(x,y)-\tilde u(x,y)$:
\begin{gather*}
Tv(x,y)=0,~~x\in\mathbb R~~ ,0<y<b,\\
v(x,0)=u(x,0)-\tilde u(x,0)=0,\\
v(x,b)=u(x,b)-\tilde u(x,b)=\\
=\frac{3}{5} b^6 x^3-\frac{1}{10} b^5 x^2-\frac{1}{20} b^9 x+\frac{1}{280} b^8
\end{gather*}
The solution of the Dirichlet problem for an inhomogeneous equation is a polynomial
\[
u(x,y)=\tilde u(x,y)+v(x,y)=
\]
\[
=\tilde u(x,y)+\frac{3b^6}{5 }v_3 (x,y)-\frac{b^5}{10} v_2 (x,y)-\frac{b^9}{20} v_1 (x,y)+\frac{b^8}{280} v_0 (x,y)=
\]
\[
=\frac{y}{840} (-504x^3 y^5+504b^5 x^3+84x^2 y^4-84b^4 x^2+
\]
\[
+42xy^8-252b^5 xy^3+210b^8 x-3y^7+14b^4 y^3-11b^7 ).
\]

\section{The Dirichlet problem for a strip in a mixed domain.}

      Consider the strip $–a<y<a$, in which the Tricomi equation has a mixed type. We give a polynomial solution of the Dirichlet problem
\[
Tu(x,y):=yu_{xx} (x,y)+u_{yy} (x,y)=0,~~-a<y<a,     \eqno{(13)}
\]
\[
u(x,-a)=\phi(x),~~   u(x,a)=\psi(x),      \eqno{(14)}
\]
where $\phi(x)$ and $\psi(x)$ are polynomials.

      In the same way as in the previous section, the solution of the Dirichlet problem (13), (14) with $\phi(x)=x^n,~\psi(x)=0$  is the function
\[
u_n (x,y)=\sum_{m=0}^{[n/2]}C_n^{2m} x^{n-2m}p_{2m}(y),
\]
where $p_{2m}(y)$ are polynomials in $y$ whose generating function is
\[
\frac{\textup{Ai}(t^{2/3} y)\textup{Bi}(t^{2/3}a)-\textup{Bi}(t^{2/3} y)\textup{Ai}(t^{2/3} a)}{ \textup{Ai}(-t^{2/3} a)\textup{Bi}(t^{2/3}a)-\textup{Bi}(-t^{2/3} a)\textup{Ai}(t^{2/3}a) }=\sum_{m=0}^{\infty}p_{2m}(y)\frac{(-1)^m t^{2m}}{(2m)!}.
\]
We give some first polynomials $p_{2m}(y)$   and the corresponding solutions of the Dirichlet problem $u_n (x,y)$. To emphasize their dependence on the parameter $a$, we will write $p_{2m}(y,a)$  
and $u_n (x,y,a)$ if necessary.
\begin{gather*}
p_0 (y)=\frac{a-y}{2a},\\
p_2 (y)=\frac{1}{12a} (y^4-2ay^3+2a^3 y-a^4 ),\\
p_4 (y)=\frac{1}{210a} (-5y^7+14ay^6-35a^3 y^4+35a^4 y^3-30a^6 y+21a^7 ),\\
p_6 (y)=\frac{1}{252a} (2y^{10}-7ay^9+30a^3 y^7-42a^4 y^6+90a^6 y^4-126a^7 y^3+103a^9 y-50a^{10} ).\\
u_0 (x,y)=\frac{a-y}{2a},\\
u_1 (x,y)=\frac{x(a-y)}{2a},\\
u_2 (x,y)=\frac{1}{12a} (-6x^2 y+6ax^2+y^4-2ay^3+2a^3 y-a^4 ),\\
u_3 (x,y)=\frac{x}{4a} (-2x^2 y+2ax^2+y^4-2ay^3+2a^3 y-a^4 ).
\end{gather*}

The solution of the Dirichlet problem (13), (14) with $\phi(x)=0,~\psi(x)=x^n$  is the function
   \[  
v_n (x,y)=\sum_{m=0}^{[n/2]}C_n^{2m} x^{n-2m}q_{2m}(y),
\]
where $q_{2m}(y)$   are polynomials in $y$ whose generating function is
\[
\frac{\textup{Ai}(t^{2/3} y)\textup{Bi}(-t^{2/3}a)-\textup{Bi}(t^{2/3} y)\textup{Ai}(-t^{2/3} a)}{ \textup{Ai}(t^{2/3} a)\textup{Bi}(-t^{2/3}a)-\textup{Bi}(t^{2/3} a)\textup{Ai}(-t^{2/3}a) }=\sum_{m=0}^{\infty}p_{2m}(y)\frac{(-1)^m t^{2m}}{(2m)!}.
\]
Since the generating function for the polynomials $q_{2m}(y)$ is obtained from the generating function for the polynomials $p_{2m}(y)$ by replacing $a$ by $-a$, the polynomials $q_{2m}(y)$ themselves are obtained from the polynomials $p_{2m} (y)$ by replacing $a$ by $-a,~ q_{2m} (y,a)=p_{2m} (y,-a)$.  The same remark applies to the corresponding solutions of the Dirichlet problem:
\[
v_n (x,y,a)=u_n (x,y,-a).
\]

\textbf{Example . } 
We consider the Dirichlet problem for the inhomogeneous Tricomy equation
\[
Tu(x,y)=2x^2 y^3-18x^3 y^4,~~x\in\mathbb R ,~~-a<y<a
\]
\[
u(x,-a)=0,~~u(x,a)=0,~~x\in\mathbb R .
\]
A particular solution of the equation is the polynomial
\[
\tilde u(x,y)=-\frac{3}{5} x^3 y^6+\frac{1}{10} x^2 y^5+\frac{1}{20} xy^9-\frac{1}{280} y^8
\]
and the Dirichlet problem for an inhomogeneous equation reduces to the Dirichlet problem for a homogeneous equation for the function $v(x,y)=u(x,y)-\tilde u(x,y)$:
\[
Tv(x,y)=0,~~x\in\mathbb R ,~~-a<y<a,
\]
\[
v(x,-a)=u(x,-a)-\tilde u(x,-a)=\frac{3}{5} a^6 x^3+\frac{1}{10} a^5 x^2+\frac{1}{20} a^9 x+\frac{1}{280} a^8,
\]
\[
v(x,a)=u(x,a)-\tilde u(x,a)=\frac{3}{5} a^6 x^3-\frac{1}{10} a^5 x^2-\frac{1}{20} a^9 x+\frac{1}{280} a^8.
\]
The solution of the Dirichlet problem for an inhomogeneous equation is a polynomial
\[
u(x,y)=\tilde u(x,y)+v(x,y)=
\]
\[
=\tilde u(x,y)+\frac{3}{5} a^6 u_3 (x,y,a)+\frac{1}{10} a^5 u_2 (x,y,a)+\frac{1}{20} a^9 u_1 (x,y,a)+\frac{1}{280} a^8 u_0 (x,y,a)+
\]
\[
+\frac{3}{5} a^6 u_3 (x,y,-a)-\frac{1}{10} a^5 u_2 (x,y,-a)-\frac{1}{20} a^9 u_1 (x,y,-a)+\frac{1}{280} a^8 u_0 (x,y,-a)=
\]
\[
=-\frac{3}{5} x^3 y^6+\frac{3}{5} a^6 x^3+\frac{1}{10} x^2 y^5-\frac{1}{10} a^4 x^2 y+\frac{1}{20} xy^9-\frac{3}{5} a^6 xy^3+\frac{11}{20} a^8 xy-\frac{1}{280} y^8+
\]
\[
+\frac{1}{60} a^4 y^4-\frac{11}{840} a^8.
\]

\section{On the uniqueness of the solution of the Dirichlet problem.}

 If we seek a solution of the Dirichlet problem for the strip $–a<y<a$ in the class of functions of polynomial growth, then the solution is not unique. To the polynomial solution indicated in section 4, we can add any solution of the form
\[
(c_1\sin(\mu_k^{3/2}x)+c_2\cos(\mu_k^{3/2}x))(\textup{Bi}(\mu_ky)\textup{Ai}(-\mu_ka)-\textup{Ai}(\mu_ky)\textup{Bi}(-\mu_ka)),
\]
where $c_1,c_2$ are arbitrary constants, and $\mu_k$ is any positive root of equation
\[
\textup{Ai}(-\mu a)\textup{Bi}(\mu a)-\textup{Ai}(\mu a)\textup{Bi}(-\mu a)=0.
\]
We give several first positive roots of this equation for $a=1$ (approximate values)
\[
\mu_1=2.340667730,~\mu_2=4.087953380,~\mu_3=5.520559835,~\mu_4=6.786708090,
\]
\[
\mu_5=7.944133587.
\]

       If we seek solutions in the class of polynomials, then this solution will be unique. To prove this it suffices to show that if the polynomial $P(x,y)$  satisfies the homogeneous Tricomi equation
\[
       TP(x,y):=yP_{xx} (x,y)+P_{yy} (x,y)=0       \eqno{(15)}
\]
and homogeneous boundary conditions
\[
                     P(x,-a)=0,~~P(x,a)=0,  \eqno{(16)}
\]
then it is identically equal to zero, $P(x,y)\equiv0$.
Substituting $P(x,y)$ in the Tricomi equation (15) and the boundary conditions (16) and equating the coefficients of monomials $x^n y^m$ to zero, we obtain that all the coefficients of the polynomial $P(x,y)$ are zero.

\section*{                                             Conclusion}
    
Using the Fourier transform of generalized functions of slow growth, we obtain polynomial solutions of the Dirichlet problem in the strip for the inhomogeneous Tricomi equation with a polynomial right-hand side and polynomial boundary conditions. By the same method, polynomial solutions can be obtained for other boundary-value problems for the Tricomi equation in the strip. For example, for the mixed boundary value problem Dirichlet-Neumann

\end{document}